\documentclass[12pt,preprint]{amsart}

\usepackage{amssymb,amsfonts,amsthm,amsmath,amscd}

\theoremstyle{plain}
\newtheorem{theorem}{Theorem}
\newtheorem*{corollary}{Corollary}
\theoremstyle{remark}
\newtheorem{remark}{Remark}

\def\j{\mathcal{J}}
\def\G{\mathcal{G}}

\begin{document}

\title{A restriction of Euclid}

\author{Grant Cairns}
\author{Nhan Bao Ho}

\address{Department of Mathematics, La Trobe University, Melbourne, Australia 3086}
\email{G.Cairns@latrobe.edu.au}
\email{nbho@students.latrobe.edu.au, honhanbao@yahoo.com}

\begin{abstract}
Euclid is a well known two-player impartial combinatorial game. A position in  Euclid is a pair of positive integers and the players move alternately by subtracting a positive integer multiple of one of the integers from the other integer without making the result negative. The player who makes the last move wins. There is a variation of  Euclid due to Grossman in which the game stops when the two entrees are equal. We examine a further variation  that we called M-Euclid in which the game stops when one of the entrees is a positive integer multiple of the other. We solve the Sprague-Grundy function for M-Euclid and compare the Sprague-Grundy functions of the three games.
\end{abstract}

\maketitle

\section{Introduction}
 \emph{Euclid} is a two-player impartial combinatorial game, introduced by Cole and Davie \cite{Cole}. In Euclid, a position is a pair of positive integers. The players move alternately, and each move is to subtract a positive integer multiple of one of the entrees from the other without making the result negative. The player who reduces one of the entrees to zero wins. In the variation of Euclid due to Grossman \cite{Gros}, the game stops when the two entrees are equal. Various aspects of Euclid and Grossman's game have been examined in the literature; see the references in \cite{GHL}.

In this note, we examine a variation, that we call \emph{M-Euclid}, in which the game stops when one of the entrees is a positive integer multiple of the other. We denote the Sprague-Grundy functions of Euclid, Grossman's game and M-Euclid  $\G_E,\G_G$ and $\G_M$ respectively. We first recall the results for $\G_E$ and $\G_G$. The convention here is that we write continued fractions $[a_0,a_1,\dots,a_n]$ so that $a_n>1$ if $n>0$.

\begin{theorem} \label{T:cf} \cite{GHL, Niv}.
Let $0<a<b$, consider the continued fraction expansion $[a_0,a_1,\dots,a_n]$ of $\frac{b}{a}$, and  let $\mathcal{I}(a,b)$ be the largest nonnegative integer $i$ such that
$a_0=\dots=a_{i-1}\leq a_{i}$.
Then
\[
\G_E(a,b)=\left\lfloor \frac{b}{a}\right\rfloor-
 \begin{cases}
0&:\  \text{if}\ \mathcal{I}(a,b)\ \text{is even},\\
 1 &:\  \text{otherwise.}\
 \end{cases}
\]
Furthermore, for Grossman's game, $\G_G(a,b)=\G_E(a,b)$ except when $a_0=a_1=\dots=a_{n}$, in which  case,
\[
\G_G(a,b)=\G_E(a,b)-(-1)^{\mathcal{I}(a,b)}.
\]
\end{theorem}


Typically, small variations in the terminal condition of a combinatorial game can produce wildly different  Sprague-Grundy functions. Interestingly, the Sprague-Grundy functions of Euclid, Grossman's game and M-Euclid are closely related. We have:

\begin{theorem}  \label{A-E T3}
Let $0<a<b$ where $b$ is not a multiple of $a$, consider the continued fraction expansion $[a_0,a_1,\dots,a_n]$ of $\frac{b}{a}$, and  let $\mathcal{J}(a,b) $ be the largest nonnegative integer $j<n$ such that $a_0=\dots=a_{j-1}\leq a_{j}$.
Then
 \[
    \G_M(a,b) = \left\lfloor \frac{b}{a} \right\rfloor -
    \begin{cases}
    0, & \text{if}\ \mathcal{J}(a,b)\ \text{is even},\\
    1, &\text{otherwise}.
    \end{cases}
    \]
\end{theorem}

\begin{remark} We draw the reader's attention to the subtle difference in the definitions of $\mathcal{I}(a,b)$ and $\mathcal{J}(a,b)$. For $\mathcal{J}(a,b)$ we have imposed $\mathcal{J}(a,b)<n$. So $\mathcal{J}(a,b)=\min\{\mathcal{I}(a,b),n-1\}$.
\end{remark}

\begin{corollary} \label{Cor}
With the notation of Theorems \ref{T:cf}  and  \ref{A-E T3},
$\G_M(a,b)=\G_E(a,b)$ except when $a_0=a_1=\dots=a_{n-1}\leq a_{n}$, in which  case,
\[
\G_M(a,b)=\G_E(a,b)-(-1)^{\mathcal{I}(a,b)}.
\]
Furthermore,
$\G_M(a,b)=\G_G(a,b)$ except when $a_0=a_1=\dots=a_{n-1}<a_{n}$, in which  case,
\[
\G_M(a,b)=\G_G(a,b)-(-1)^{\mathcal{I}(a,b)}.
\]
\end{corollary}

Having found the right formulation of  Theorem \ref{A-E T3}, its proof is straight-forward. We follow closely the proof of  \cite[Theorem 1]{GHL}.

This paper continues our investigations of variations of Euclid and related questions; see \cite{GH1,GH2,GH3,GHL}.

\section{Proof of Theorem \ref{A-E T3}}

For convenience we write $\G$ instead of $\G_M$ and by abuse of language, we write $\j(p)$ and $\G(p)$ for their values at a position $p=[a_0,a_1,\dots,a_n]$.  It suffices to establish the following two  properties:
\begin{enumerate}
\item For every move $p\mapsto q$, we have $\G(q)\not=\G(p)$.
\item If $\G(p)>0$, then for all integers $k$ with $0\leq k< \G(p)$, there exists a move $p\mapsto q$ such that $\G(q)=k$.
 \end{enumerate}

We will make repeated use of the following fact: if $p=[a_0,a_1,\dots,a_n]$ and $\j(p)$ is odd, then $a_0\leq a_1$ and $n>1$; indeed, if $n=1$ or $a_0> a_1$, then we would have  $\j(p)=0$. Similarly, if $\j(p)$ is even then either $a_0\geq a_1$ or $n=1$.

First observe that  Theorem \ref{A-E T3} holds for $n=1$. Indeed, clearly $\G(1,a_1) = 1$ for all $a_1$ and hence by induction, $\G([a_0,a_1]) = a_0$ for all $n$. Since we have $a_0= \lfloor \frac{b}{a} \rfloor$, and $\mathcal{J}(a,b)=0$, the result follows. So we need only deal with positions $p$ having $n>1$.

To establish (1), suppose we have a move $p\mapsto q$ with $\G(q)=\G(p)$. First suppose that $q=[a_0-i,a_1,\dots,a_n]$ for some $1\leq i<a_0$.  From the definition of $\G$, it is clear that  $i=1$, $\j(p)$ is odd and $\j(q)$ is even. As $\j(p)$ is odd, $a_0\leq a_1$, and so as $\j(q)$ is even, $a_0-1\geq a_1$. Hence $a_0\leq a_1\leq a_0-1$, which is impossible. So we may assume that $q=[a_1,\dots,a_n]$. At first sight, as $\G(q)=\G(p)$, there are three possibilities:
\begin{enumerate}
\item[(i)] $a_0=a_1-1$ and $\j(p)$ is even and $\j(q)$ is odd,
\item[(ii)] $a_0=a_1+1$ and $\j(p)$ is odd and $\j(q)$ is even,
\item[(iii)] $a_0=a_1$ and $\j(p)$ and $\j(q)$ have the same parity.
\end{enumerate}
But case (i) is impossible, since $a_0\geq a_1$ when $\j(p)$ is even, case (ii)  is impossible since $a_0\leq a_1$ when $\j(p)$ is odd, and case (iii) contradicts the definition of $\j$.

To establish (2), suppose that $0\leq k< \G(p)$.   First suppose  that $\j(p)$ is odd, so $\G(p)=a_0-1$.  Consider the position $q=[k+1,a_1,\dots,a_n]$. Since $\j(p)$ is odd, $a_0\leq a_1$.  In particular,  $k+1<a_1$ and thus $\j(q )=1$. It follows that $\G(q)=k$, as required. So it remains to treat the case where  $\j(p)$ is even.  In this case, $\G(p)=a_0$ and $a_0 \geq a_1$.

We first treat the situation where $k=0$. Assume for the moment that $a_0>1$. Consider  $q=[1,a_1,\dots,a_n]$. Notice that we may assume that $\j(q)$ is even, since otherwise $\G(q)=0$, as required.  In particular, we have  $a_1=1$. Let $q'=[a_1,\dots,a_n]$. But if $\j(q)$ is even, then $\j(q')$ is odd and hence $\G(q')=a_1-1=0$, as required. Similarly, if $a_0=1$, then as $\j(p)$ is even, we have  $a_1=1$, and since $\j(p)$ is even, $\j(q')$ is odd and $\G(q')=0$. This
completes the case $k=0$.

Now suppose that $0< k< \G(p)$ and let $q=[k,a_1,\dots,a_n]$. If $\j(q)$ is even, then $\G(q)=k$, as required. So we may assume that $\j(q)$ is odd and thus $k\leq a_1$.  In this case, we have $\G(q)=k-1$.  Let $q'=[k+1,a_1,\dots,a_n]$.  If $\j(q')$ is odd, then $\G(q')=k$, as required, so we may assume that $\j(q')$ is even, and therefore $k+1\geq a_1$.  Thus $k+1\geq a_1\geq k$. Hence, either $k+1= a_1$ or  $k= a_1$. Consider $q''=[a_1,\dots,a_n]$. If $k=a_1$, then as $\j(q)$ is odd,   $\j(q'')$ is even, and hence $\G(q'')=a_1=k$, as required.
Finally, if $k+1=a_1$,  then as $\j(q')$ is even, $\j(q'')$ is odd, and hence $\G(q'')=a_1-1=k$, as required.


\end{document}